\documentclass[12pt,a4paper]{article}
\usepackage{t1enc}
\usepackage[latin1]{inputenc}
\usepackage[english]{babel}
\usepackage{amsmath}
\usepackage{amsfonts}
\usepackage{graphicx}
\usepackage{mathptmx}

\font\bBB=msbm10

\def\bBR{\mbox{\bBB R}}

\begin{document}
\title {A Note on Topologically-Trivial Braids}
\author{
 Orlin Stoytchev\thanks{ American University in Bulgaria, 2700 Blagoevgrad, Bulgaria; \newline\indent Institute for Nuclear Research and Nuclear Energy, 1784 Sofia , Bulgaria; \newline\indent ostoytchev@aubg.bg} }
\date{}
\maketitle
\abstract{We give a simple characterization of braids that can be unplaited keeping separately their upper ends and their lower ends tied together}
\section*{}
Consider Artin's \cite{Artin} braid group $B_n$ and its subgroup of pure braids $P_n\subset B_n$. Each representative of $P_n$ can be regarded as a geometric object embedded in $\bBR^3$ in the following way: two {\bf pieces} of horizontal planes, one being a vertical translate of the other, and $n$ strands each connecting a point on the lower plane with the translate of that point on the upper plane. The strands of course do not intersect each other and in addition the vertical coordinate of each strand is a monotonic function of the parameter of the strand. Now consider those elements of $P_n$ that are topologically trivial as embedded in $\bBR^3$, namely there are isotopies bringing them to the trivial braid. These are not just the usual isotopies for Artin braids where the two parallel planes are fixed in $\bBR^3$. We allow moves like rotations of the planes, flipping a strand above (below) and around one of the planes, passing the plane through some part of the braid. The elements of $P_n$ satisfying this property obviously form a subgroup $R_n\subset P_n$. This subgroup was studied by Shepperd \cite{Shep} using algebraic methods. By bringing each pure braid to a certain canonical form, the author provides a (rather complicated) algorithm of deciding whether or not this braid can be unplaited keeping its ends tied together, i.e., whether or not it is topologically trivial in the above sense.\par
In the present note we exhibit a simple and intuitive solution of that same problem. It can be described in words as follows: Pick one of the strands, e.g., the $n$-th strand, and pass the whole braid, starting with the upper common end, through itself along the $n$-th strand, then straighten it  back up again. You get a (pure) braid in which the $n$-th strand runs straight, without crossing any of the other strands. Then the original braid is topologically trivial if and only if by
removing the $n$-th strand we get an  $n-1$ braid belonging to the center
of $P_{n-1}$ or, equivalently, represents some power of the full twist (see below) of  $n-1$ strands.\par
Shepperd \cite{Shep} considered the following moves for plaiting a braid, while keeping the ends tied together: Pass the common (upper) end of the braid between the first $k$ strands and the remaining $n-k$ strands and pull the common end back up again. Depending on the direction of rotation of the common end, the two groups of strands will be twisted each around itself in opposite directions. Let us denote by $b_k$ the element of $P_n$ obtained in this way from the trivial braid. In terms of the Artin's generators
$\sigma_i$ we have:
\begin{equation}
b_k=(\sigma_{k-1}\cdots \sigma_2\sigma_1  )^{-k} \, ( \sigma_{n-1} \cdots\sigma_{k+2}\sigma_{k+1})^{n-k}\ .
\end{equation}
The element $b_0\equiv d$ is a full (clockwise) twist of the whole braid, the element $b_1$ is a braid in which the first strand is straight, while the remaining $n-1$ strands are twisted around themselves, the element $b_n\equiv d^{-1}\equiv b_0^{-1}$ is a full twist in the opposite (counterclockwise) direction.
One may consider a seemingly more general move in which one picks arbitrarily $k$ strands from the braid, pulls them to the left, then passes the common end between them and the remaining $n-k$ strands. However this more general move is equivalent to first applying a suitable element of $B_n$ which rearranges the strands, then splitting 
the strands into the first $k$ and the remaining $n-k$, then passing the common end between these groups and finally applying the inverse element of $B_n$. As we shall see, the subgroup generated by $\{b_k\}$ is normal in $B_n$, so it is enough to apply the simpler moves, producing the elements $b_k$. For the same reason
it is enough to apply these operations only to the upper end of the braid --- in other words you do not get anything more by passing the common end
through a nontrivial braid somewhere in the middle. Thus the question of finding the topologically trivial braids reduces to the study of the subgroup $R_n$, generated by the elements $\{b_k\}$.\par\noindent
{\bf Note.} The statement that there are no more general moves than the ones described above, allowing one to unplait a braid while keeping its ends tied together, should be taken as a (very plausible) conjecture. It is such also in the work of Shepperd \cite{Shep}. \par
\begin{figure} [h]
\setlength{\abovecaptionskip}{20pt}
\setlength{\belowcaptionskip}{0pt}
\centering
\includegraphics[scale=0.28]{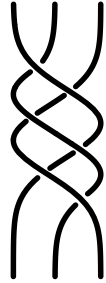}
\caption{The full twist $d$ in the case $n=3$ }
\label{twist}
\end{figure}
\begin{figure} [h]
\setlength{\abovecaptionskip}{20pt}
\setlength{\belowcaptionskip}{0pt}
\centering
\includegraphics[scale=0.45]{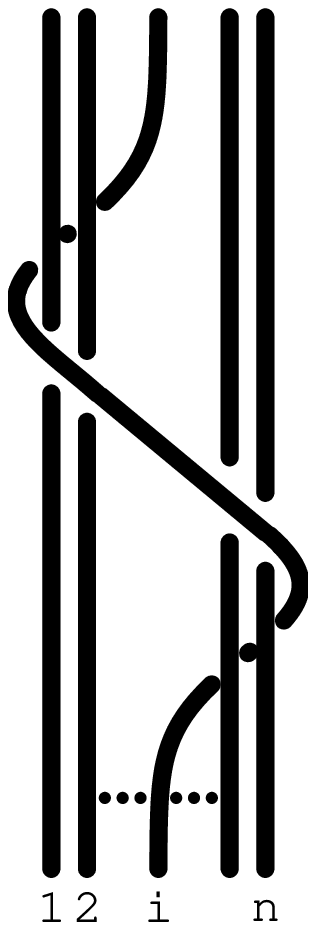}
\caption{The flip $r_i$ }
\label{flip}
\end{figure}

For our purposes it is more convenient to use a different set of elements, generating the same subgroup $R_n$, namely
\begin{equation}{\begin{array}{l}
d:=(\sigma_{n-1}\cdots \sigma_2\sigma_1)^n,\\
r_1:=\sigma_1\sigma_2\cdots \sigma_{n-2}\sigma_{n-1}^2\sigma_{n-2}\cdots \sigma_1,\\
r_2:=\sigma_1^2\sigma_2\cdots\sigma_{n-2}\sigma_{n-1}^2\sigma_{n-2}\cdots \sigma_2,\\
r_i:=\sigma_{i-1}\cdots\sigma_2\sigma_1^2\sigma_2\cdots\sigma_{n-2}\sigma_{n-1}^2\sigma_{n-2}\cdots\sigma_i,
\quad i=2,3,\ldots n-1,\\
r_n:=\sigma_{n-1}\sigma_{n-2}\cdots\sigma_2\sigma_1^2\sigma_2\cdots\sigma_{n-2}\sigma_{n-1}.\\
\end{array}}\label{gener}
\end{equation}
We call $d$ a {\it full twist} (Fig.\ref{twist}) and we call $r_i$ {\it flips} (Fig.\ref{flip}).
It is a simple matter to see that the following identity relates the generators $b_k$ to the generators $r_i$ and $d$:
\begin{equation}
b_k=d(r_k r_{k-1}\cdots r_1)^{-1}=(r_k r_{k-1}\cdots r_1)^{-1} d\ \ .
\end{equation}
Notice that the move described by $r_k r_{k-1}\cdots r_1$ is similar to the one corresponding to $r_1$ but performed to the whole bunch of the first $k$ strands. Using algebra (Artin's braid relations) one sees that a flip of all $n$ strands gives the same effect as two full twists:
\begin{equation}
(r_n r_{n-1}\cdots r_2 r_1)=d^2=( r_{n-1}r_{n-2}\cdots  r_1 r_n)=\cdots =(r_1 r_n\cdots r_3 r_2) \ \ .
\label{dsquare}
\end{equation}
The last equation is a mathematical expression of the so-called ``belt trick'' demonstrating the well-known fact in topology that a complete rotation by $720^{\circ}$ is homotopic to the identity. (An explanation of the simple connection between the fundamental group of $SO(3)$ and braids can be found in \cite{Sto}.) 
Equation \ref{dsquare} shows that the subgroup $R_n\subset P_n$ is actually generated by $\{r_i\}_{i=1,\ldots n-1}$ and $d$. We have the following
\par\bigskip\noindent
{\bf Lemma 1}\par\noindent
{\it The subgroup $R_n$ is normal in $B_n$.}\par\smallskip\noindent
{\bf Proof:} Since $d$ is central in $B_n$ it suffices to exhibit explicit formulas for the conjugates of all flips $r_i$. The following identities can be checked directly:
\begin{equation}{\begin{array}{l }
\sigma_j r_i \sigma_j^{-1}=\sigma_j^{-1}r_i \sigma_j=r_i,\quad i-j>1\ \hbox{or}\ j-i>0,\\
\sigma_{i-1}r_i\sigma_{i-1}^{-1}=r_i r_{i-1} r_i^{-1},\\
\sigma_{i-1}^{-1} r_i \sigma_{i-1}=r_{i-1},\\
\sigma_i r_i \sigma_i^{-1}=r_{i+1},\quad i\le n-1\\
\sigma_i^{-1} r_i \sigma_i = r_i^{-1}r_{i+1} r_i,\quad i\le n-1\\
\end{array}}\label{normal}
\end{equation}\hfill $\Box$
\par\smallskip
We will denote by $R_n'$ the subgroup generated by $r_i,\ i=1,\dots , n$. (This is obviously normal in $B_n$ and contains $d^2$ but not $d$.) 
Let us observe that the $n$ elements $r_i,\ i=1,\dots n$ satisfy the cyclicity relations contained in Equation \ref{dsquare}, but if we limit ourselves to the (first) $n-1$ of the $r_i$, they are probably free.
The next proposition formalizes the process, described in the second paragraph --- passing a braid through itself along the $n$-th strand:\par\bigskip\noindent
{\bf Proposition 1} \par\noindent
{\it Each class in the factorgroup $P_n/R_n'$ contains a representative in which the $n$-th strand runs straight, without being entangled with the others. This representative is unique up to a full twist of the first $n-1$ strands.}\par\smallskip\noindent
{\bf Proof:} To prove the existence we describe algebraically a process (see examples below) through which from each braid $b\in P_n$ we get a new braid $s(b)$ with the prescribed property by multiplying $b$ with elements from $R_n'$. Let $b$ be written as some word in the Artin generators $\sigma_i,\ \ i=1, \dots , n-1$ and their inverses:
\begin{equation}
b=\sigma_{i_1}^{\pm 1}\sigma_{i_2}^{\pm 1}\cdots \sigma_{i_k}^{\pm 1}\ ,\quad i_j\in \{1,2,\dots n-1\}\ \ .
\end{equation}
Now we track the original $n$-th strand, starting from the right-hand end of the word, and mark those letters, which correspond to this strand passing in front of another strand. There are two possibilities --- the original $n$-th strand is at $i$-th position and we come across a letter $\sigma_{i-1}$ (the current $i$-th strand passes in front of the $(i-1)$-st) or we come across a letter $\sigma_i^{-1}$ (the current $i$-th strand passes in front of the $(i+1)$-st). In the first case we insert on the left of the marked letter an element $r_i^{-1}$, in the second case we insert on the left of the marked letter an element $r_i$. We call the new braid obtained after all insertions $s(b)$. In this braid the original $n$-th strand passes always behind the other strands, so the word can be simplified to represent a pure braid in which the $n$-th strand runs straight, or, equivalently, the corresponding word does not contain $\sigma_{n-1}$.\par
In order to prove uniqueness and for the further development we need the concept of a {\it spherical braid} --- several distinct points on a sphere and the same number of points, in the same positions, on a smaller sphere, connected by strands in such a way that the radial coordinate of each strand is monotonic in the parameter of that strand.
When one considers the obvious multiplication for spherical braids and  isotopy classes one obtains the so called {\it braid group of the sphere} \cite{Fad2}, which algebraically is $B_n/\{r_i\}_{i=1,\dots n}$. This is known as the mapping-class group of the sphere (with 0 punctures and 0 boundaries). Let's denote this by $S_n$. Using a stereographic projection one can map an Artin braid to a spherical braid and vice-versa (choosing an axis for the projection not intersecting any strand). In the forward direction (from Artin braids to spherical ones) this induces a homomorphism between the braid groups:
$$
\pi:B_n\rightarrow S_n
$$
and we have precisely 
$$
Ker\, \pi=R_n'\ .
$$
The moves described by the generators $r_i$ (see Eq. \ref{gener}) for the corresponding spherical braid represent flipping the $i$-th strand around the inner sphere (which is an isotopy for the spherical braid) while $d$ is a rotation of the inner sphere by $360^{\circ}$ around some axis, e.g., the one we used to define the stereographic projection). Once we have obtained a representative in which the $n$-th strand runs straight (as in the existence part of the proof) we see that the only moves of the corresponding spherical braid, leaving this strand straight, are rotations by $360^{\circ}$ around this strand and isotopies in the usual (Artin) sense of the $n-1$ strands. \hfill $\Box$
\par\bigskip
Applying topological arguments we see that the following is true:
\par\smallskip\noindent
{\bf Proposition 2} \par\noindent
{\it A braid $b\in P_n$ is topologically trivial (i.e., can be unplaited keeping separately its lower and upper ends tied together) if and only if the corresponding topologically equivalent braid $s(b)$ constitutes a full twist of the first $n-1$ strands. The topologically nonequivalent braids are in 1-1 correspondence with the elements of $P_{n-1}/\{d\}$, i.e., the elements of the pure braid group of $n-1$ strands modulo its center.}
\par\smallskip\noindent
{\bf Proof:} We consider our braid as a spherical braid. It is topologically trivial if and only if it can be brought to the trivial spherical braid by isotopies of the strands and possibly full rotations around arbitrary axes. In fact, taking into account the topology of $SO(3)$ and its connections with spherical braids (see \cite{Sto}) it is enough to take a single rotation by $360^{\circ}$ around some fixed axis (or no rotation at all). Next we observe that if there is an isotopy bringing a  spherical braid to the trivial one then there is another isotopy which first straightens the $n$-th strand and then, keeping it straight brings the remaining $n-1$ strands to a configuration where they are untangled except possibly twisted around each other. This rather intuitive fact relies on choosing a couple of continuous maps $\sigma : [0,1]\times [0,1]\rightarrow SO(3)$ satisfying certain boundary conditions. To avoid too many notations we prefer to skip the details. Note that even though there may be an isotopy bringing directly the initial braid to the trivial one, by first straightening the $n$-th strand we may end up with a braid for which the remaining $n-1$ strands are twisted some even number of times around each other. Finally, we can untwist those strands by a rotation around the $n$-th strand. \hfill $\Box$
\par\bigskip\noindent
{\bf Example 1 -- the Standard three-strand braid:}
\par\noindent
The standard three-strand braid as well as any three-strand braid can be unplaited with its ends tied together, as is well-known. This is a simple corollary of Proposition 2 since $P_2$ --- the pure braid group of 2 strands --- is the cyclic group generated by the corresponding full twist $d$.
\par\bigskip\bigskip\noindent
\begin{figure} [h]
\setlength{\abovecaptionskip}{20pt}
\setlength{\belowcaptionskip}{0pt}
\centering
\includegraphics[scale=0.25]{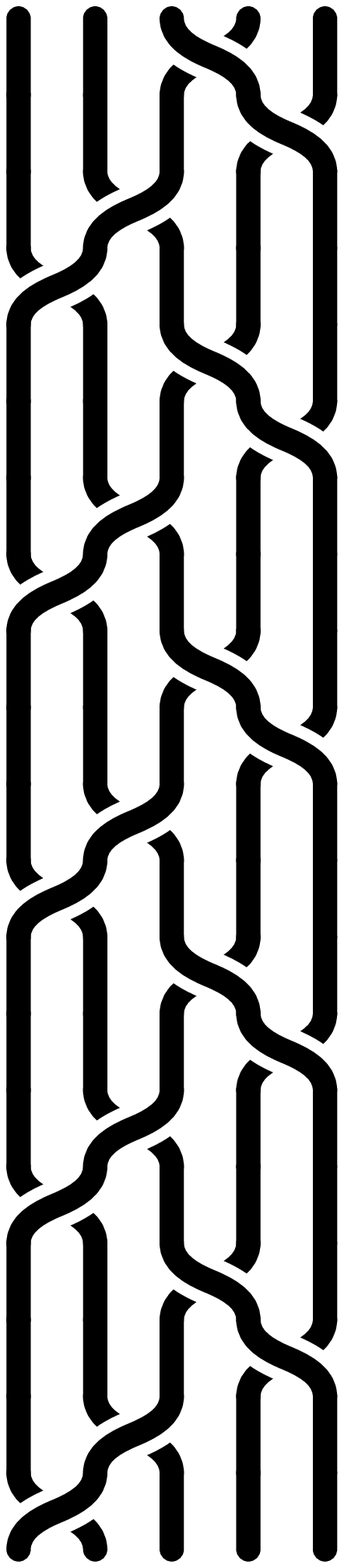}\quad\quad\quad\quad \quad\quad \includegraphics[scale=0.21]{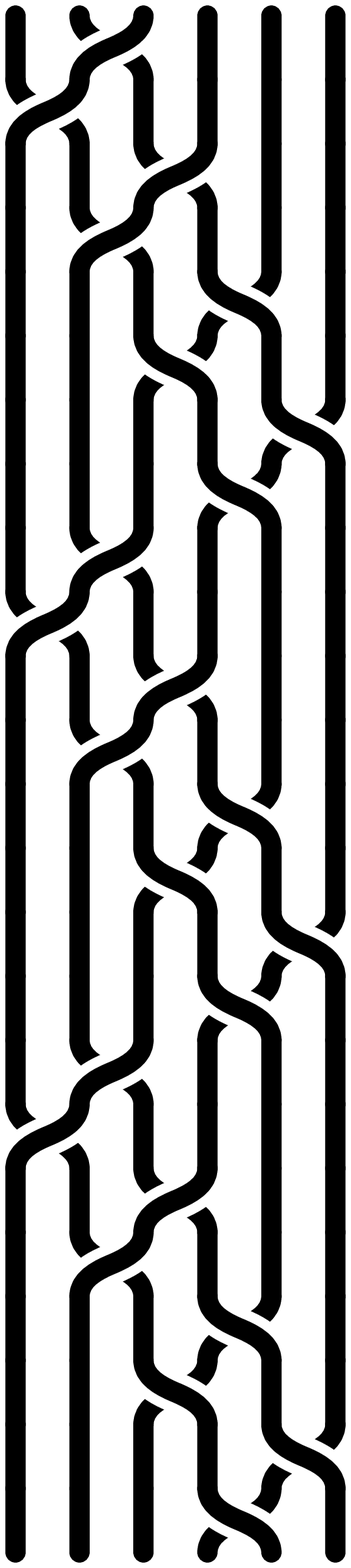}
\caption{(a) The English sennit\quad\quad (b) The braided theta }
\label{esennit}
\end{figure}
{\bf Example 2 -- the English sennit:}\par\noindent
 Consider the braid with $5$ strands given by the following word (Fig.\ref{esennit}(a)):
\begin{equation}
b=(\sigma_3\sigma_4\sigma_2^{-1}\sigma_1^{-1})^5\ \ .
\end{equation}
We expand the power and put asterisks above those letters, corresponding to the fifth strand passing in front of another strand:
\begin{equation}
b=\sigma_3\sigma_4\sigma_2^{-1}\sigma_1^{-1}\       \sigma_3\sigma_4\overset{*}{\sigma}_2^{-1}\overset{*}{\sigma}_1^{-1}\    \sigma_3\sigma_4\sigma_2^{-1}\sigma_1^{-1}\   \sigma_3\sigma_4\sigma_2^{-1}\sigma_1^{-1}\   \overset{*}{\sigma}_3\overset{*}{\sigma}_4\sigma_2^{-1}\sigma_1^{-1}\ \ .
\end{equation}
The transformed braid $s(b)$ is obtained by inserting appropriate elements $r_i^{\pm 1}$ to the left of the marked letters:
\begin{eqnarray}
s(b)&=&\sigma_3\sigma_4\sigma_2^{-1}\sigma_1^{-1}      \sigma_3\sigma_4r_2{\sigma}_2^{-1}r_1{\sigma}_1^{-1}    \sigma_3\sigma_4\sigma_2^{-1}\sigma_1^{-1}  \sigma_3\sigma_4\sigma_2^{-1}\sigma_1^{-1}  r_4^{-1}{\sigma}_3r_5^{-1}{\sigma}_4\sigma_2^{-1}\sigma_1^{-1}\nonumber\\ 
&=&\sigma_3\sigma_4\sigma_2^{-1}\sigma_3\sigma_4\sigma_1\sigma_2\sigma_3\sigma_4^2\sigma_3\sigma_1\sigma_2\sigma_3\sigma_4^2\sigma_3
\sigma_2\sigma_3\sigma_4\sigma_2^{-1}\sigma_1^{-1}\sigma_3\sigma_4\sigma_2^{-1}\sigma_1^{-1}\times \nonumber \\
&\times& \sigma_4^{-2}\sigma_3^{-1}\sigma_2^{-1}\sigma_1^{-2}\sigma_2^{-1}\sigma_4^{-1}\sigma_3^{-1}\sigma_2^{-1}\sigma_1^{-2}\sigma_2^{-1}
\sigma_3^{-1}\sigma_2^{-1}\sigma_1^{-1} \ .\nonumber
\end{eqnarray}
In principle one can simplify directly the above word and show that it represents the trivial braid. It is somewhat easier to remove the fifth strand and work with the remaining braid with $4$ strands. Indeed, tracking the fifth strand in $s(b)$ we see that it crosses all other strands always behind, so the braid $s(b)$ is equivalent to
the braid obtained by removing the fifth strand and then adding it running straight at fifth position. The braid obtained by removing the fifth strand from $s(b)$ can be described as follows: Suppose that the fifth strand is currently at $j$-th position and you encounter a letter $\sigma_i^{\pm 1}$ with $i\ne j, j+1$ (in other words this generator does not move the current $j$-th strand). Define a new letter
$$
m_j(\sigma_i^{\pm 1})=\begin{cases}
\sigma_i^{\pm 1}\quad \text{if}\ \  i<j\\ 
\sigma_{i-1}^{\pm 1}\quad  \text{if}\ \  i>j+1
\end{cases}
$$
and extend this transformation multiplicatively to a whole part of the word for which the original fifth strand does not change its position. In the original word $s(b)$ cross out those letters that correspond to a move of the original fifth strand. The rest of the word splits into pieces for which that strand does not move. Apply the transformation
$m_j$ with the appropriate $j$ to each piece. In our example we obtain the following braid:
\begin{eqnarray}
s'(b)
&=&\sigma_3\mathaccent\times{\sigma}_4\sigma_2^{-1}\mathaccent\times{\sigma}_3 m_3(\sigma_4\sigma_1)\mathaccent\times{\sigma}_2
m_2(\sigma_3\sigma_4^2\sigma_3)\mathaccent\times{\sigma}_1m_1(\sigma_2\sigma_3\sigma_4^2\sigma_3
\sigma_2\sigma_3\sigma_4\sigma_2^{-1})\mathaccent\times{\sigma}_1^{-1}m_2(\sigma_3\sigma_4)\times\nonumber \\
&\times & \mathaccent\times{\sigma}_2^{-1}m_3(\sigma_1^{-1}
\sigma_4^{-2})\mathaccent\times{\sigma}_3^{-1}\sigma_2^{-1}\sigma_1^{-2}\sigma_2^{-1}
\mathaccent\times{\sigma}_4^{-1}\sigma_3^{-1}\sigma_2^{-1}\sigma_1^{-2}\sigma_2^{-1}
\sigma_3^{-1}\sigma_2^{-1}\sigma_1^{-1} \nonumber\\
&=&\sigma_3\sigma_2^{-1}\sigma_3\sigma_1\sigma_2\sigma_3^2\sigma_2\sigma_1\sigma_2\sigma_3^2\sigma_2\sigma_1\sigma_2\sigma_3
\sigma_1^{-1}\sigma_2\sigma_3\sigma_1^{-1}
 \sigma_3^{-2}\sigma_2^{-1}\sigma_1^{-2}\sigma_2^{-1}\times\nonumber\\
&\times&\sigma_3^{-1}\sigma_2^{-1}\sigma_1^{-2}
\sigma_2^{-1}\sigma_3^{-1}\sigma_2^{-1}\sigma_1^{-1}\ \ .\nonumber
\end{eqnarray}
Now, using Artin's relations after somewhat lengthy calculations we get
$$
s'(b)=Id\ ,
$$
which means that $b\in R_5$ or in other words the English sennit is topologically trivial.\par\bigskip\noindent
{\bf Example 3 -- the braided theta:} 
\par\noindent
Consider the braid with $6$ strands given by the following word (Fig. \ref{esennit}(b)):
\begin{equation}
b=(\sigma_2^{-1}\sigma_1^{-1}\sigma_3^{-1}\sigma_2^{-1}\sigma_4\sigma_3\sigma_5\sigma_4)^3
\end{equation}
This is like the normal three-strand braid except that each strand in the latter is replaced by two strands running together. Therefore we can think of each such pair as a ribbon. Thus we have something like the classical three-strand braid but plaited out of ribbons, which run flat without being twisted. You can take a strip of paper, cut two longitudinal slits to form a ``theta'' and then try to plait it to obtain the ``braided theta'' (see Bar-Natan's gallery of knotted objects \cite{Natan} from which the example was borrowed). We have
\[
b=\sigma_2^{-1}\sigma_1^{-1}\sigma_3^{-1}\sigma_2^{-1}\sigma_4\sigma_3\sigma_5\sigma_4
\sigma_2^{-1}\sigma_1^{-1}\overset{*}{\sigma}_3^{-1}\overset{*}{\sigma}_2^{-1}\sigma_4\sigma_3\sigma_5\sigma_4
\sigma_2^{-1}\sigma_1^{-1}\sigma_3^{-1}\sigma_2^{-1}\overset{*}{\sigma}_4\sigma_3\overset{*}{\sigma}_5\sigma_4\ , 
\]
\begin{eqnarray}
s(b)&=&\sigma_2^{-1}\sigma_1^{-1}\sigma_3^{-1}\sigma_2^{-1}\sigma_4\sigma_3\sigma_5\sigma_4
\sigma_2^{-1}\sigma_1^{-1}r_3\sigma_3^{-1}r_2\sigma_2^{-1}\sigma_4\sigma_3\sigma_5\sigma_4
\sigma_2^{-1}\sigma_1^{-1}\sigma_3^{-1}\sigma_2^{-1}r_5^{-1}\sigma_4\sigma_3 r_6^{-1}\sigma_5\sigma_4 \hfil\nonumber \\
&=&\sigma_2^{-1}\sigma_1^{-1}\sigma_3^{-1}\sigma_2^{-1}\sigma_4\sigma_3\sigma_5\sigma_4\sigma_2^{-1}\sigma_1^{-1}\sigma_2\sigma_1^2
\sigma_2\sigma_3\sigma_4\sigma_5^2\sigma_4\sigma_1^2\sigma_2\sigma_3\sigma_4\sigma_5^2\sigma_4\sigma_3\times\hfil\nonumber\\
&\times& \sigma_4\sigma_3\sigma_5\sigma_4\sigma_2^{-1}\sigma_1^{-1}\sigma_3^{-1}\sigma_2^{-1}\sigma_5^{-2}\sigma_4^{-1}
\sigma_3^{-1}\sigma_2^{-1}\sigma_1^{-2}\sigma_2^{-1}\sigma_5^{-1}\sigma_4^{-1}\sigma_3^{-1}\sigma_2^{-1}\sigma_1^{-2}\sigma_2^{-1}\sigma_3^{-1}\ .
\hfil\nonumber
\end{eqnarray}
Removing the sixth strand produces the following:
\begin{eqnarray}
s'(b)&=&\sigma_2^{-1}\sigma_1^{-1}\sigma_3^{-1}\sigma_2^{-1}\sigma_4\sigma_3\mathaccent\times{\sigma}_5\mathaccent\times{\sigma}_4
m_4(\sigma_2^{-1}\sigma_1^{-1}\sigma_2\sigma_1^2
\sigma_2)\mathaccent\times{\sigma}_3m_3(\sigma_4\sigma_5^2\sigma_4\sigma_1^2)\mathaccent\times{\sigma}_2
m_2(\sigma_3\sigma_4\sigma_5^2\sigma_4\sigma_3 \sigma_4\times\hfil\nonumber\\
&\times&\sigma_3\sigma_5\sigma_4)\mathaccent\times{\sigma}_2^{-1}m_3(\sigma_1^{-1})\mathaccent\times{\sigma}_3^{-1}
m_4(\sigma_2^{-1}\sigma_5^{-2})\mathaccent\times{\sigma}_4^{-1}
\sigma_3^{-1}\sigma_2^{-1}\sigma_1^{-2}\sigma_2^{-1}\mathaccent\times{\sigma}_5^{-1}\sigma_4^{-1}\sigma_3^{-1}\sigma_2^{-1}\sigma_1^{-2}\sigma_2^{-1}\sigma_3^{-1}
\hfil\nonumber\\
&=&\sigma_2^{-1}\sigma_1^{-1}\sigma_3^{-1}\sigma_2^{-1}\sigma_4\sigma_3\sigma_2^{-1}\sigma_1^{-1}\sigma_2\sigma_1^2
\sigma_2\sigma_3\sigma_4^2\sigma_3\sigma_1^2\sigma_2\sigma_3\sigma_4^2\sigma_3\sigma_2\sigma_3\sigma_2
\sigma_4\sigma_3\times\hfil\nonumber\\
&\times&
\sigma_1^{-1}\sigma_2^{-1}\sigma_4^{-2}\sigma_3^{-1}\sigma_2^{-1}\sigma_1^{-2}\sigma_2^{-1}\sigma_4^{-1}\sigma_3^{-1}\sigma_2^{-1}\sigma_1^{-2}\sigma_2^{-1}\sigma_3^{-1}
\hfil\nonumber\\
&=&Id
\end{eqnarray}
This shows that what we have called ``braided theta'' is topologically trivial as a braid of $6$ strands with the ends tied together. In order to show that the ``braided theta'' is trivial when considered as plaited out of three ribbons we need more --- we need to check that all transformations involved move the strands in pairs. This is indeed the case. Alternatively we can work directly with the ribbons (pairs of strands) and introduce  {\it ribbon flips} $R_1$, $R_2$, $R_3$ and their inverses, which are similar to the ones in Figure \ref{flip} but performed on the 3 ribbons.
By definition we have $R_i:=r_{2i}  r_{2i-1}$ and the effect of a flip $R_i$ is similar to that of the usual flip $r_i$ except that it twists the $i$-th ribbon by 
$720^{\circ}$ (counterclockwise). It is easier to find experimentally, rather than doing the algebra, that the ``braided theta'' in Figure \ref{esennit} (b) is the product
$R_3 R_2^{-1} R_3^{-1} R_2$.
\par
 
\vfill\eject
\end{document}